\documentclass{amsart}
\usepackage{amssymb}
\usepackage{adjustbox}
\usepackage[margin=1.2in]{geometry}

\usepackage{epsfig}  		% For postscript
\usepackage{epic,eepic}       % For epic and eepic output from xfig
\usepackage{tikz}

\newtheorem{theorem}{Theorem}[section]
\newtheorem{proposition}[theorem]{Proposition}
\newtheorem{lemma}[theorem]{Lemma}

\theoremstyle{definition}

\theoremstyle{remark}
\newtheorem{remark}[theorem]{Remark}

\numberwithin{equation}{section}

\begin{document}

\title{A Fixed--Point Approach to Non--Commutative Central Limit Theorems}
\author{Jad Hamdan}
\address{Mathematical Institute, University of Oxford}
\email{hamdan@maths.ox.ac.uk}

%\author{Second Author}
%\address{Department of Mathematics, University of South Carolina,
%Columbia, SC 29208}
%\email{second@math.sc.edu}
%\urladdr{www.math.sc.edu/$\sim$second}

\begin{abstract}
We show how the renormalization group approach can be used to prove quantitative central limit theorems (CLTs) in the setting of free, Boolean, bi--free and bi--Boolean independence under finite third moment assumptions. The proofs rely on the construction of a contraction on a subspace of probability measures over $\mathbb{R}$ (or $\mathbb{R}^2$) equipped with a suitable metric, which has the appropriate analogue of a Gaussian distribution as a fixed point (for instance, the semi--circle law in the case of free independence).  In all cases, this yields a convergence rate of $1/\sqrt{n}$, and we show that this can be improved to $1/n$ in some instances under stronger assumptions.
\end{abstract}

\maketitle

\section{Introduction}

In non--commutative probability, one deals with elements of a $*$--algebra which do not necessarily commute. This gives rise to many distinct notions of independence, and in turn, to different binary operations on probability measures, in the same way as classical independence gives rise to convolution. The classical theory of sums of independent random variables is then, quite often, perfectly paralleled by these new theories, which have their own counterparts of the central limit theorem (CLT), L{\'e}vy--Khintchine formula and other well--known results. 

A prototypical example of this is the theory of free independence (leading to free additive convolution), which was first introduced by Voiculescu in \cite{voic} and has since been heavily studied, culminating in an explicit correspondence between limit laws for classical and free additive convolution established in a seminal work of Bercovici and Pata \cite{BercoviciPata}. Another well--understood and much simpler example is that of Boolean convolution, introduced by Speicher and Woroudi in \cite{SpeicherWoroudi}. In a series of works \cite{Voiculescu1,Voiculescu2}, Voiculescu also introduced an extension of free probability for pairs of algebras, which enables the study of non--commutative left and right actions of algebras on a reduced free product space, simultaneously. This theory of so--called \textit{bi--freeness} has since attracted much attention and been rapidly developed (see \cite{CharlesworthNelsonSkoufranis1,CharlesworthNelsonSkoufranis2,GuHuangMingo, HuangWang, SkoufranisFree1, SkoufranisFree2}), including a recently--established Bercovici--Pata--type bijection \cite{HasebeHuangWang}. The theory of Boolean convolution was similarly extended to pairs of unital algebras in \cite{GuSkoufranis}.

The problem of quantifying the rate of convergence in these various central limit theorems (CLTs) (as accomplished by Berry and Esseen in \cite{Berry, Esseen} in the classical setting) has also been the subject of much interest. A convergence rate for the free CLT was first obtained by Kargin \cite{Kargin, Kargin2} for compactly supported measures, and this assumption was later weakened by Chistyakov and G{\"o}tze \cite{ChistyakovGotze1, ChistyakovGotze2}. Berry--Esseen--type results for Boolean convolution were first established by Arizmendi and Salazar \cite{ArizmendiSalazar} and refined by the latter in \cite{Salazar}. Rates were also obtained in multidimensional and operator--valued settings \cite{SpeicherMultiFree, SpeicherMultiBE, Banna, MaiSpeicher, ArizmendiBannaTseng}. All of these results are in terms of the Kolmogorov or L{\'e}vy metric. 

Following an expository work of Ott \cite{ott} which presents a renormalization group proof of the classical CLT (originally due to Hamedani \cite{hamedani}), this paper shows how this type of argument can be adapted to non--commutative settings. In particular, we prove quantitative versions of the free, Boolean, bi--free and bi--Boolean central limit theorems with a decay of $1/\sqrt{n}$ in all cases (see theorems \ref{FreeBE}, \ref{BoolBE}, \ref{BifBE}, and \ref{BboolBE} respectively). These are the first Berry--Esseen--type results for the latter two. 

\bigskip

The main idea behind our arguments will be to define a metric $d$ on a subspace of the space of probability measures on $\mathbb{R}$ (or $\mathbb{R}^2$) for which a {renormalization} map $T$ fixes the appropriate analogue of the Gaussian distribution, and is a contraction. Our definition of $d$ here is analogous to that of the Fourier--based metric used in Ott's proof \cite{ott} (of which an in--depth discussion can also be found in \cite{jose}), except that we replace the Fourier transform in its definition by the $R$--transform in the free setting, and the self--energy in the Boolean setting. 

The map $T$ that we use to prove the free CLT was first used by Anshelevich in  \cite{Anshelevich1998TheLO}, where it is referred to as the \textit{central limit operator}. To the best of our knowledge, Anshelevich's is the only previous work to explicitly interpret the free CLT as a fixed point theorem, but their approach and results differ significantly from ours. Namely, instead of introducing a metric, they focus on the problem of linearizing $T$ (computing its G{\^a}teaux derivative). This allows them to compute its eigenvalues, and show that its eigenfunctions are absolutely continuous with respect to the Lebesgue measure, with densities equalling multiples of Chebyshev polynomials of the first kind. The central limit theorem then follows as a corollary: in a subspace of compactly supported measures, the differential of $T$ evaluated at $\rho_{\text{sc}}$ has spectrum inside the unit disk, and we then expect $\rho_{\text{sc}}$ to be an attracting fixed point on that subspace. By contrast, the proof below applies to a much larger family of measures.

The paper is organized as follows. Sections 2 and 4 give the necessary background and definitions to study the free/Boolean and bi-free/bi--Boolean central limit theorems, respectively, from an analytic viewpoint. In particular, we introduce the (single and double--variable) $R$--transform, Cauchy transform and self--energy, describe their asymptotic behavior near zero and give characterizations of weak convergence in terms of these analytic functions. Using this, we prove our main results in sections 3 and 5.
\section{Preliminary results and notation}

Let $\mathcal{M}=\mathcal{M}(\mathbb{R})$ denote the set of (Borel) probability measures on $\mathbb{R}$. We denote the $k$--th moment of a measure $\mu$ by $m_k(\mu)$, and let $\mathcal{M}^k:=\{\mu\in \mathcal{M}\,:\, |m_j(\mu)|<\infty \, \forall j\leq k\}$. For every $k\geq 1$, let
\[
    \mathcal{M}_0^k := \left\{\mu \in \mathcal{M}^k \,:\, m_1(\mu)=0,\,\,m_2(\mu)=1\right\}.
\]
% and use the superscript $\mathcal{M}^c$ to denote the measures in $\mathcal{M}$ which are compactly supported. 

For any $\mu,\nu\in\mathcal{M}$, we denote by $\mu\boxplus \nu$ the \textit{free additive convolution} and by $\mu\uplus\nu$ the \textit{Boolean convolution} of the measures $\mu$ and $\nu$. These are the distributions of $X+Y$ when $X$ and $Y$ are free/Boolean independent random variables with laws $\mu$ and $\nu$ respectively, but we will opt for purely analytic definitions which are given later.

Let $\mathcal{NC}(k)$ and $\mathcal{I}(k)$ be the lattices of non--crossing partitions and interval partitions of $\{1,...,k\}$ respectively. For measures with compact support (and for which all moments are finite), we define the \textit{free cumulants} $(\kappa_n(\mu))_{n\geq 1}$ (resp. \textit{Boolean cumulants} $(r_n(\mu))_{n\geq 1}$) of $\mu$ by
\begin{equation}\label{momentcumulantfree}
    m_k(\mu) = \sum_{\substack{\pi \in \mathcal{NC}(k), \\\pi = \{B_1,...,B_n\}} } \prod_{j=1}^n \kappa_{|B_j|}(\mu),\quad \bigg(\text{resp.} \quad   m_k(\mu) = \sum_{\substack{\pi \in \mathcal{I}(k), \\\pi = \{B_1,...,B_n\} }} \prod_{j=1}^n r_{|B_j|}(\mu)\bigg).
\end{equation}

Note the similarity between these formulae and the classical moment--cumulant formulae, where the sum is taken over all partitions. Much like the latter, these relations can be inverted using the theory of M{\"o}bius functions of lattices (see the book of Stanley \cite{stanley}). One can show that free cumulants linearize the free convolution of measures, satisfying
\begin{equation}\label{cumulantconv}
     \kappa_n(\mu\boxplus \nu) = \kappa_n(\mu)+\kappa_n(\nu) \quad (n\geq 1).
\end{equation}
The same is true for {Boolean convolution} and its associated cumulants. In fact, one can define free/Boolean convolution to be the composition law for which such a linearization holds, but these definitions would not generalize to unbounded measures. This problem is circumvented by translating the above moment--cumulant relations into the analytic statements of the next section. 

\subsection{Analytic theory of free and Boolean convolution.}\label{analytictheory}

The \textit{Cauchy transform} of a probability measure $\mu\in\mathcal{M}(\mathbb{R})$ is the analytic function on the upper half--plane (minus the support of $\mu$)
\[
    G_\mu:\,z\mapsto\int_\mathbb{R}d\mu(t)/(z-t).
\]
As a result of the \textit{Stieltjes inversion formula}, every probability measure on $\mathbb{R}$ is uniquely determined by its Cauchy transform. The $R-$transform $R_{\mu}$ of $\mu$ is then defined as the analytic solution to 
\begin{equation}\label{Rtransf}
     G_\mu\Big(R_\mu(z)+\frac{1}{z}\Big)=z
\end{equation}
whose domain of definition will depend on the assumptions made on $\mu$. If none are made, this will be the union of a truncated \textit{Stolz angle }$\Delta_{\alpha,\beta}\subseteq \mathbb{C}^{-}$ {(where $\mathbb{C}^{-}$ denotes the lower half of $\mathbb{C}$)}, defined by
\[
    \Delta_{\alpha,\beta} = \{x+iy\in\mathbb{C}^{-}\,:\, |x|< -\alpha y, y>-\beta\}
\]
(for some $\alpha,\beta>0$) with its complex conjugate $\overline{\Delta_{\alpha,\beta}}:=\{\bar{z}\,:\,z\in \Delta_{\alpha,\beta}\}$.
We define the \textit{self--energy} of $\mu$ to be the function $E_\mu(z)=z-1/G_\mu(z)$. Note that both $E_\mu$ and $R_\mu$ uniquely determine $\mu$ since they can be used to recover $G_\mu$. Lastly, since $G_\mu(\bar{z})=\overline{G_\mu(z)}$, one can easily check that these properties are inherited by $E_\mu$ and $R_\mu$. 

\begin{remark}
A well--defined inverse for $G$ (and in turn, an $R$--transform) exists for all probability measures, following the work of Voiculescu \cite{VoiculescuRTransform} in the compactly supported case, and its subsequent generalization by Maassen \cite{Maassen}, Chistyakov and G{\"o}tze \cite{ChistyakovGotze1,ChistyakovGotze2}, and others \cite{BercoviciVoiculescu,BelinschiBercovici}. 
A key step in this generalization is the use of subordination functions (see, for instance, Chapter 3 of \cite{MingoSpeicher}). 
\end{remark}

For any two probability measures $\mu_1$ and $\mu_2$ on $\mathbb{R}$ with $R$--transforms $R_1$ and $R_2$, there exists a unique probability  measure $\nu$ whose $R$-transform is $R_1+R_2$ (see \cite{MingoSpeicher} for a proof). We define $\mu_1\boxplus \mu_2$ to be this measure. Following \cite{SpeicherWoroudi} we similarly define $\mu_1\uplus \mu_2$ to be the measure whose self--energy is $E_{\mu_1\uplus \mu_2}(z)=E_{\mu_1}(z)+E_{\mu_2}(z)$.

The more we assume about the probability measure $\mu$, the better behaved its $R$--transform is. In particular, if $\mu$ is compactly supported with support in the interval $[-r,r]$, then $R_\mu$ is analytic in a disc centered around $0$ with radius $1/(6r)$. Moreover, the coefficients in this expansion are the aforementioned free cumulants  $(\kappa_n(\mu))_{n\geq 1}$ of $\mu$, making $R_\mu$ their generating function. Similarly, the coefficients in the expansion of $E_\mu(z)$ for large enough $z$ are the Boolean cumulants, but we note that the latter is in negative powers of $z$ and thus view $E_\mu(1/z)$ as the natural Boolean analogue of $R_\mu$. 

 If $\mu$ isn't compactly supported but has finite second moment $\sigma^2$, then $R_\mu(z)$ is analytic on a disc with center $-i/(4\sigma)$ and radius $1/(4\sigma)$. Since $0$ is on the boundary of this disc, we may not have free cumulants beyond the second. However, if $\mu$ has a moment of order $p$, a result of Benaych--Georges (Theorem 1.3 in \cite{BenaychGeorges}) gives a Taylor expansion with $p$ terms. 
 \begin{theorem}[Benaych-Georges]\label{BGtheorem} Let $p$ be a positive integer and $\mu$ a probability measure on the real line. If $\mu$ admits a $p$--th moment, then $R_\mu$ admits the Taylor expansion
 \[
    R_\mu(z)=\sum_{i=0}^{p-1}\kappa_{i+1}(\mu)z^i+o(z^{p-1})
 \]
 where $(\kappa_n(\mu))_{n\in\mathbb{N}}$ are the free cumulants of $\mu$ \eqref{momentcumulantfree} and the limit is as $z\to 0$ non--tangentially, meaning $|z|\to 0$ and $|\Re(z)|\leq -\alpha \Im{z}$ for some $\alpha>0$.
 \end{theorem}

The analogous result for the self--energy is the following (proposition 13 in \cite{SalazarThesis}).
\begin{proposition}[Arizmendi--Salazar]\label{Kexpansion} Let $p$ be a positive integer and $\mu$ a probability measure on the real line. If $\mu$ admits a $p$--th moment, then $E_\mu$ admits the expansion
\[
    E_\mu(z)=\sum_{i=0}^{p-1}\frac{r_{i+1}(\mu)}{z^i}+o\bigg(\frac{1}{z^{p-1}}\bigg) 
\]
where $(r_n(\mu))_{n\in\mathbb{N}}$ are the Boolean cumulants of $\mu$ and the limit is as $z\to \infty$ non--tangentially.
\end{proposition}

We end this section with the two observations, beginning with the following scaling property which follows directly from the definitions of $R_\mu$ and $E_\mu$ and is true for any $\lambda \in \mathbb{R}$. With the convention $(\lambda \mu)(A):=\mu(A/\lambda)$ for $A\subseteq \mathbb{R}$,
\begin{equation}\label{scalingRK} 
    R_{\lambda \mu}(z)=\lambda R_{\mu}(\lambda z), \quad E_{\lambda \mu}(1/z)=\lambda E_{\mu}(1/\lambda z)
\end{equation}
Let $\rho_{sc}, \rho_b$ denote the semi-circle and Bernoulli distributions, defined as
\[
    \frac{1}{2\pi}\sqrt{\text{max}(4-x^2,0)}\mathrm{d}x\,\text{ and } \,\frac{1}{2}(\delta_{1}+\delta_{-1})
\]
respectively (where $\mathrm{d}x$ is the Lebesgue measure on $\mathbb{R}$). Then using (\ref{scalingRK}) and equating cumulants (which are $\kappa_n(\rho_{sc})=r_n(\rho_b)=\mathbf{1}_{n=2}$), we obtain that 
 \begin{equation}\label{stability}
    R_{(\rho_{sc}\boxplus\rho_{sc})/\sqrt{2}}(z)=R_{\rho_{sc}}(z)=z, \quad E_{(\rho_b\uplus\rho_b)/\sqrt{2}}(1/z)=E_{\rho_b}(1/z)=z.
\end{equation}

\subsection{Weak convergence via analytic transforms} We will need the following characterizations of weak convergence in terms of the $R$-transform and self-energy. Both propositions are immediate corollaries of their two dimensional counterparts in section \ref{bisection} (theorems \ref{charbi} and \ref{charbibool}).
\begin{proposition}\label{Rconv}
        Let $\{\nu_n\}_{n\geq 1}$ be a sequence of Borel probability measures. Then $\nu_n$ converges weakly to a probability measure $\nu$ on $\mathbb{R}$ if and only if 
    \begin{enumerate}
        \item There exists a Stolz angle $\Delta$ such that all $R_{\nu_n}$ are defined on $\Delta$.
        \item $\lim_{n\to\infty}R_{\nu_n}(z)=R_\nu(z)$ for every $z\in \Delta$, and $R_{\nu_n}(-iy)\to 0$ uniformly in $n$ as $y\to 0^+$. 
    \end{enumerate}
\end{proposition}
\begin{proposition}\label{Econv} Let $\{\nu_n\}_{n\geq 1}$ be a sequence of Borel probability measures.  Then $\nu_n$ converges weakly to a probability measure $\nu$ on $\mathbb{R}$ if and only if 
 $\lim_{n\to \infty}E_{\nu_n}(z)=E_\nu(z)$ for all $z\in(\mathbb{C}\setminus\mathbb{R})$, and the limit $E_{\nu_n}(z)\to 0$ holds uniformly in $n$ as $|z|\to\infty$ non--tangentially.
\end{proposition}

\section{Fixed-point proofs of the free and Boolean CLTs}
Fix $\epsilon>0$ small enough so that $\Delta_{\epsilon,\epsilon}\subseteq \{|z+i/4|<1/4\}$. For any $k\geq 2$, define the distances,
\begin{align}\label{metrics}
    d^{(k)}_{\text{Free}}(\mu,\nu)=\sup_{z\in \Delta_{\epsilon,\epsilon}}\frac{|R_\mu(z)-R_\nu(z)|}{|z|^k},\quad  d^{(k)}_{\text{Bool}}(\mu,\nu)=\sup_{z\in \Delta_{\epsilon,\epsilon}}\frac{|E_\mu(1/z)-E_\nu(1/z)|}{|z|^k}.\end{align}
Throughout this section, we omit the superscript from $d^{(k)}_{\text{Free}}$ and $d^{(k)}_{\text{Bool}}$ when $k=2$. We will prove that $d_{\text{Free}}, d_{\text{Bool}}$ are finite metrics on $\mathcal{M}_0^{3}(\mathbb{R})$, and that convergence with respect to these metrics implies weak convergence. The following theorems will then follow straightforwardly.
\begin{theorem}[Free Berry--Esseen]\label{FreeBE} Let $\mu\in \mathcal{M}_0^3$. Then there exists a $C(\mu)=C>0$ such that
\[
    d_{\text{\normalfont Free}}(\tfrac{1}{\sqrt{n}}\mu^{\boxplus n}, \rho_{sc}) \leq \frac{|m_3(\mu)|+C}{\sqrt{n}}.
\]
In particular, $\tfrac{1}{\sqrt{n}}\mu^{\boxplus n}$ converges weakly to $\rho_{sc}$. 
If one assumes that $m_3(\mu)=0$ and $m_4(\mu)<\infty$, this can be improved to $({|m_4(\mu)|+C})/n$ for the distance $d^{(3)}_{\text{\normalfont{Free}}}$ and a different constant $C>0$.
\end{theorem}
\begin{theorem}[Boolean Berry--Esseen]\label{BoolBE} Let $\mu\in \mathcal{M}_0^3$. Then there exists a $C(\mu)=C>0$ such that
\[
    d_{\text{\normalfont{Bool}}}(\tfrac{1}{\sqrt{n}}\mu^{\uplus n}, \rho_b) \leq \frac{|m_3(\mu)|+C}{\sqrt{n}}.
\]
In particular, $\tfrac{1}{\sqrt{n}}\mu^{\uplus n}$ converges weakly to $\rho_{b}$.
If one assumes that $m_3(\mu)=0$ and $m_4(\mu)<\infty$, this can be improved to $({|m_4(\mu)|+C)/n}$ for the distance $d^{(3)}_{\text{\normalfont{Bool}}}$ and a different constant $C>0$.
\end{theorem}
\begin{remark}
    A rate of $1/\sqrt{n}$ in Kolmogorov distance was obtained by Chistyakov and G{\"o}tze \cite{ChistyakovGotze1} for the free CLT assuming a finite {absolute} third moment. The same authors then showed that the rate could be improved to $1/n$ if $m_3(\mu)=0, m_4(\mu)<\infty$, a fact which is mirrored by the theorem above. For the Boolean CLT, a rate of $1/\sqrt{n}$ (in L{\'e}vy distance) has recently been obtained by Salazar \cite{Salazar} for measures with finite sixth moment. Assuming a finite fourth moment, an earlier work of Arizmendi and Salazar \cite{ArizmendiSalazar} obtained a rate of $1/n^{1/3}$.
\end{remark}
Define the \textit{renormalization map} with respect to free (resp. Boolean) convolution to be the map $T^{\boxplus}: \mu\mapsto (\mu\boxplus \mu)/\sqrt{2}$ (resp. $T^{\uplus}: \mu\mapsto (\mu\uplus \mu)/\sqrt{2}$). 
Then by (\ref{scalingRK}), $$R_{T^{\boxplus}\mu}(z)=\sqrt{2}R_\mu(z/\sqrt{2}),\quad E_{T^{\uplus}\mu}(1/z)=\sqrt{2}E_\mu(1/(\sqrt{2}z)).$$
Analyzing the first few coefficients in the partial Taylor expansions reveals that $\mathcal{M}_0^3$ is closed under the action of $T^\boxplus$ and $T^\uplus$, and (\ref{stability}) can be rewritten as $ T^{\boxplus}\rho_{sc} =\rho_{sc}, T^{\uplus}\rho_b = \rho_b.$

Theorems \ref{FreeBE} and \ref{BoolBE} will follow straightforwardly from the following propositions.

\begin{proposition}\label{qe:finite1d}
    $d_{\text{\normalfont Free}}$ and $d_{\text{\normalfont Bool}}$ are finite metrics on $\mathcal{M}_0^3$, where convergence in the metric topology implies weak convergence.
    In particular, for any $\mu\in \mathcal{M}_0^3$, there exists constants $B, C>0$ such that
    \begin{align}\label{bound}
         d_{\text{\normalfont Free}}(\mu, \rho_{sc})\leq |m_3(\mu)|+B,\quad d_{\text{\normalfont Bool}}(\mu, \rho_b)\leq |m_3(\mu)|+C.
    \end{align}
    If $\mu\in \mathcal{M}_0^{4}$ and $m_3(\mu)=0$, the claim holds for the $d^{(3)}$ distances as well, replacing $m_3(\mu)$ by $m_4(\mu)$ in the right hand side of the inequalities.
\end{proposition}
\begin{proof} 
    For both $d_{\text{\normalfont Free}}$ and $d_{\text{\normalfont Bool}}$, symmetry is clear, and separation follows from the identity theorem (using complex conjugation to extend to the whole of $\mathbb{C}\setminus \mathbb{R}$) and the fact that probability measures are uniquely determined by their Cauchy/$R$--transform. The triangle inequality follows from $\sup f+g \leq \sup f + \sup g$ and the triangle inequality for the complex norm.

    The fact that convergence in $d_{\text{Free}}$ implies weak convergence is an immediate consequence of theorem \ref{Rconv}. For $d_{\text{Bool}}$, we must also argue that $d_{\text{Bool}}(\mu_n,\mu)\to 0$ implies $E_{\mu_n}\to E_\mu$ on all of $(\mathbb{C}\setminus \mathbb{R})$, for $\mu_n,\mu\in \mathcal{M}_0^3$. To that end, note that for any $R>1$, Chebyshev's inequality gives $\mu_n([-R,R])\geq 1-1/R^2:=C>0$ uniformly over $n$. For any compact $K\subseteq \mathbb{C}\setminus \mathbb{R}$ and $z\in K$, it follows that
    \[
        |\Im G_{\mu_n}(z)|\geq |y|\int_{|t|\leq R} \frac{1}{(x-t)^2+y^2}\mathrm{d}\mu_{n}(t)\geq C \min_{z=x+iy\in K} \frac{|y|}{(|x|+R)^2+y^2} >0
    \]
     uniformly over $n\geq 1$ and $z\in K$. By Montel's theorem $\{E_{\mu_n}\}_n$ is therefore a normal family: on any compact $K\subseteq (C\setminus \mathbb{R})$, every subsequence of $E_{\mu_n}$ has a further converging subsequence to some limit, which must coincide with $E_{\mu}$ by the identity theorem and $d(\mu_n,\mu)\to 0$.

    For finiteness of $d_{\text{Free}}$ on $\mathcal{M}_0^3$, note that since the $R$--transform of a probability measure with unit variance is analytic on $|z+i/4|<1/4$, it is bounded on any compact subset of $\Delta_{\epsilon,\epsilon}$ away from zero. It therefore suffices to show that $\lim_{z\to 0} {|R_\mu(z)-R_\nu(z)|}{|z|^{-2}}<\infty$
    for any $\mu,\nu\in\mathcal{M}_0^3$ (with the limit taken inside $\Delta_{\epsilon,\epsilon}$) but this follows immediately from theorem \ref{BGtheorem}. 
    To get (\ref{bound}), we use the same theorem to write $R_\mu(z)=z+m_3(\mu)z^2+z^2v(z)$ for some $v$ satisfying $|v(z)|\to 0$ as $|z|\to 0$, and take $B=\sup_{\Delta_{\epsilon,\epsilon}}|v|<\infty$. Finiteness of $d_{\text{Bool}}$ and the upper bound for $d_{\text{Bool}}(\mu,\rho_b)$ follow by an identical argument, using the expansion for $E_\mu$ (proposition \ref{Kexpansion}).
    
    Lastly, taking one additional term in the expansions of $R_\mu$ and $E_\mu$, the proofs of the claims for $d^{(3)}_{\text{\normalfont{Free}}}$ and $d^{(3)}_{\text{\normalfont{Bool}}}$ are identical.
\end{proof}
% \begin{remark}
%     The constants $B$ and $C$ can be made arbitrarily small if one modifies the definition of the distances, taking the supremum over $y\in (0,\epsilon]$ for small $\epsilon$.
% \end{remark}

\begin{proposition}\label{contract}
    $T^{\boxplus}$ (resp. $T^{\uplus}$) is a contraction on $(\mathcal{M}_0^3,d_{\text{\normalfont Free}})$ (resp. $(\mathcal{M}_0^{3},d_{\text{\normalfont Bool}})$) with contraction constant $ 2^{-1/2}$, and 
    \begin{align}\label{iter}
        d_{\text{\normalfont Free}}(T^{\boxplus n}\nu, \rho_{sc})\leq2^{-n/2}d_\text{\normalfont Free}(\nu,\rho_{sc}), \\
        d_{\text{\normalfont Bool}}(T^{\uplus n}\nu, \rho_b)\leq2^{-n/2}d_\text{\normalfont Bool}(\nu,\rho_b)\label{iter2}.
    \end{align}
    On the subspace $\mathcal{M}_0^4\cap \{m_3(\mu)=0\}$ (equipped with the appropriate metric), the contraction constant can be improved to $1/2$ in both cases.
\end{proposition}

\begin{proof} We prove this in the free case, the Boolean case being essentially identical. Recall that the $R-$transform of $T^{\boxplus}\nu$ is $R_{T^{\boxplus}\nu}(z)=\sqrt{2}R_\nu(z/\sqrt{2})$, thus
    \begin{align*}
        d_{\text{Free}}(T\nu, T\mu)&= \sup_{z\in \Delta_{\epsilon,\epsilon}}\frac{\sqrt{2}}{2}\frac{|R_\mu(z/\sqrt{2})-R_\nu(z/\sqrt{2})|}{|z/\sqrt{2}|^2}= \sup_{z\in {2}^{-1/2}\Delta_{\epsilon,\epsilon}}\frac{1}{\sqrt{2}}\frac{|R_\mu(z)-R_\nu(z)|}{|z|^2}\leq\frac{d_{\text{Free}}(\nu, \mu)}{\sqrt{2}}.
    \end{align*}
    For $d^{(3)}_{\text{Free}}$, the increased exponent in the denominator incurs an additional factor of $1/\sqrt{2}$.
    The inequalities (\ref{iter}) and (\ref{iter2}) then follow from the fact that $T^{\boxplus}$ fixes $\rho_{sc}$.
\end{proof}
%  \begin{remark}
%     If $\mu,\nu$ and $\gamma$ are probability measures, $*$ denotes (classical) convolution and $\lambda>0$, then a metric $d$ satisfying $d(\lambda\mu, \lambda\nu)= \lambda^s d(\mu,\nu)$ and $d(\mu * \gamma, \nu*\gamma)\leq d(\mu,\nu)$ is said to be \textit{s--ideal} \cite{zolotarev}. 
%     More generally, if we let $d_s$ be the generalization of $d$ with $|z|^s$ in the denominator, then $d_s$ is $(s+1)$--ideal. This is in contrast with the analogous Fourier metric (see \cite{ott}) which is only $s$--ideal, as we are aided here by the additional factor of $\lambda$ incurred by lemma \ref{scalarR}. Indeed, \ref{contract} still holds with $s=1+\epsilon$ instead of $2$, but we chose the latter for aesthetic reasons.
% \end{remark}
\begin{remark}
	Nothing is said here about the completeness of these metric spaces, which is not needed for our main argument. To the author's best knowledge, such Cauchy/$R$-transform-based metrics have yet to be studied, by contrast with the family of Fourier--based metrics in \cite{ott} which are rather well--understood (see \cite{jose}). 
\end{remark}
The next and final proposition gives a few useful properties for our metrics. 
\begin{proposition}\label{properties}
    Let $\mu,\nu,\gamma$ and $\eta$ be measures on $\mathbb{R}$ with mean zero and variance $\leq 1$. Then 
    \begin{align}
    d_{\text{\normalfont Free}}^{(k)}(\lambda \mu, \lambda\nu) &\leq \lambda^{k+1}d_{\text{\normalfont Free}}^{(k)}(\mu, \nu),\label{propscale}
\end{align}
for any $\lambda\in (0,1)$, and whenever $m_2(\mu\boxplus\nu),m_2(\gamma\boxplus\eta)\leq 1$,
\begin{align}
    d_{\text{\normalfont Free}}^{(k)}(\mu\boxplus \nu, \eta \boxplus \gamma)&\leq d_{\text{\normalfont Free}}^{(k)}(\mu,\eta)+d_{\text{\normalfont Free}}^{(k)}(\nu, \gamma). \label{propadd}
\end{align}
Note that in both of these inequalities, the left and right-hand sides may be equal to infinity.
The analogous properties hold for and $d^{(k)}_{\text{\normalfont Bool}}$.
\end{proposition}
\begin{proof}
    (\ref{propadd}) follows from the triangle inequality (for $|\cdot |$) and definition of the $R$--transform/self--energy, while (\ref{propscale}) is an immediate consequence of equation (\ref{scalingRK}).
\end{proof}
\begin{remark}
    Metrics satisfying (\ref{propadd}) and for which (\ref{propscale}) is an equality are referred to as \textit{$3$--ideal} (more generally, $s$--ideal where $s$ is the exponent of $\lambda$ on the right hand side). The existence of such metrics was originally shown by Zolotarev \cite{zolotarev1,zolotarev}. 
\end{remark}

\subsection{Proof of theorems \ref{FreeBE} and \ref{BoolBE}} \label{proof1d}
   Let $\mu\in\mathcal{M}_0^3$. We begin with the free CLT (theorem \ref{FreeBE}), noting that its proof is immediate along geometric subsequences $N=2^n$, since
   \[
        d_{\text{Free}}(\tfrac{1}{\sqrt{N}}\mu^{\boxplus N}, \rho_{sc})=d_{\text{Free}}(T^{\boxplus N}\mu, T^{\boxplus N} \rho_{sc}) \leq N^{-1/2} d_{\text{Free}}(\mu,\rho_{sc})
   \]
   by proposition \ref{contract} and $d_{\text{Free}}(\mu,\rho_{sc})\leq |m_3(\mu)|+C$ by proposition \ref{qe:finite1d}.  To extend to more general subsequences, we use proposition \ref{properties} and the fact that $\frac{1}{\sqrt{n}}\rho_{sc}^{\boxplus n}=\rho_{sc}$ to write
\[
    d_{\text{\normalfont Free}}(\tfrac{1}{\sqrt{n}}\mu^{\boxplus n}, \rho_{sc}) = d_{\text{\normalfont Free}}(\tfrac{1}{\sqrt{n}}\mu^{\boxplus n}, \tfrac{1}{\sqrt{n}}\rho_{sc}^{\boxplus n})\leq \frac{n}{n^{3/2}}{d_{\text{\normalfont Free}}(\mu, \rho_{sc})}.\]
The proof for the Boolean case is identical, replacing every occurrence of $\boxplus$ with $\uplus$, $d_{\text{Free}}$ with $d_{\text{Bool}}$ and $\rho_{sc}$ with $\rho_b$. The improved rate of convergence in $\mathcal{M}_0^4\cap \{m_3(\mu)=0\}$ follows from using $d^{(3)}_{\text{Free}}/d_{\text{Bool}}^{(3)}$.

\section{Bi--free harmonic analysis}\label{bisection}

 Following \cite{Voiculescu1, Voiculescu2}, there is a ``two--faced" extension of free probability that enables the study of non--commutative left and right actions of algebras on a reduced free product space simultaneously. This gives rise to the notion of \textit{bi--free} independence for pairs of non--commutative random variables (which reduces to freeness when one restricts one side to be constant), and in turn to a new type of convolution on measures on $\mathbb{R}^2$. Once again, this \textit{bi--free additive convolution $\mu {\boxplus\boxplus} \nu$} is linearized by a set of cumulants relying on so--called bi--non--crossing partitions, but can be defined more generally by purely analytic means (see, e.g., \cite{Voiculescu1} and \cite{CharlesworthNelsonSkoufranis2} for combinatorial developments of the theory, and \cite{HuangWang} for their analytic counterparts). Much like in the free case, this theory has been shown to mirror the classical theory, complete with a theory of bi--free infinite divisibility \cite{HasebeHuangWang, GuHuangMingo}. Likewise, the theory of Boolean convolution was also generalized to pairs of unital algebras in \cite{GuSkoufranis}.

As a final illustrative example, we prove Berry--Esseen--type results in these bi--free and bi--Boolean settings, beginning by introduction the requisite analytic machinery.

\subsection{Bi-free $R$ and Cauchy transforms}

We first extend the definition the Cauchy transform to include Borel planar probability measures $\mu$. Let $\mathcal{M}(\mathbb{R}^2)$ be the space of such measures, 
    $m_{k,l}(\mu) = {\iint_{\mathbb{R}^2}} x^ky^l d\mu(x,y)$
be the mixed moments of $\mu$ for $k,l\geq 0$,  and
\[
    \mathcal{M}_{0,c}(\mathbb{R}^2) := \big\{\mu\in\mathcal{M}(\mathbb{R}^2)\,:\, m_1(\mu^{(i)})=0,  m_2(\mu^{(i)})=1 \mbox{ for $i=1,2,$ }\, m_{1,1}(\mu)=c\big\}.
\]
Letting $\mathcal{M}^{3}(\mathbb{R}^2)$ be the space of measures $\mu\in \mathcal{M}(\mathbb{R}^2)$ for which $\iint_{k+\ell\leq 3}|x|^k|y|^\ell d\mu<\infty$ if $k+l\leq 3$, we set $\mathcal{M}_{0,c}^3(\mathbb{R}^2):=\mathcal{M}^3(\mathbb{R}^2)\cap \mathcal{M}_{0,c}(\mathbb{R}^2)$. Note that unlike in the previous section, we require boundedness of \textit{absolute} moments here. The Cauchy transform of $\mu $ is the analytic function
\[  
    G_\mu(z,w) = \int_{\mathbb{R}^2}\frac{d\mu(s,t)}{(z-s)(w-t)}
\]
 on $(\mathbb{C}\setminus \mathbb{R})^2$, which we note satisfies $G_\mu(\bar{z},\bar{w})=\overline{G_\mu(z,w)}$. As in the single variable case, one can recover the underlying measure by Stieltjes inversion. The \textit{(bi--free partial) $R$--transform} of $\mu$ is then defined as 
\begin{equation}\label{biR}
    R_\mu(z,w)=1+zR_{\mu^{(1)}}(z)+wR_{\mu^{(2)}}(w)-\frac{zw}{G_\mu(R_{\mu^{(1)}}(z)+1/z,R_{\mu^{(2)}}(w)+1/w)},
\end{equation}
and also uniquely determines $\mu$ (see proposition 2.5 in \cite{HuangWang}).
For this to be well--defined at $(z,w)$, one must ensure that the $R$--transforms of the marginal distributions are defined at this point and that the denominator of the rightmost term never vanishes. We know that the former is true on $\Delta\cup\overline{\Delta}$ for some Stolz angle $\Delta$ depending on $\mu$. As for the nonvanishing of the Cauchy transform, we have the following asymptotic behaviour (see \cite{HuangWang}) 
\[
    G_\mu(z,w)=\frac{1}{zw}(1+o(1)) \mbox{ as $z, w\to \infty$ non--tangentially}.
\]
Since $1/\lambda+R_{\mu^{(j)}}(\lambda)=(1/\lambda)(1+o(1))$ for $j=1,2$, one can thus shrink $\Delta$ if need be to make $R_\mu$ well--defined on some product domain $\Omega=(\Delta\cup\overline{\Delta})\times (\Delta\cup\overline{\Delta})$, on which it will be holomorphic. The \textit{partial self--energy} of $\mu$ is defined by 
\[
    E_{\mu}(z,w)=\frac{1}{z}E_{\mu^{(1)}}(z)+\frac{1}{w}E_{\mu^{(2)}}(w)+\frac{G_\mu(z,w)}{G_{\mu^{(1)}}(z)G_{\mu^{(2)}}(w)}-1
\]
and is considerably simpler than its free counterpart, being defined on the entirety of $(\mathbb{C}\setminus \mathbb{R})^2$. 

If $\mu$ is compactly supported, then $E_\mu$ and $R_\mu$ admit an absolutely convergent bivariate power series expansion around $(0,0)$, with real coefficients which are the bi--free/bi--Boolean cumulants (which we do not define here). For our purposes, we will only need the following partial expansions.
\begin{proposition}[Voiculescu]\label{biRexpansion}Let $|c|\leq 1$ and $\mu \in \mathcal{M}_{0,c}^{3}(\mathbb{R}^2)$, then there exist coefficients $\{\kappa_{k,l}\}_{k+l=3}$ such that
\[
    R_\mu(z,w)=z^2+w^2+czw+\bigg(\sum_{\substack{k+l=3\\k,l \geq 0}}\kappa_{k,l}z^kw^l + o(z^kw^l)\bigg)
\]
as $|z|,|w|\to 0$ non--tangentially.
\end{proposition}

\begin{proposition}[Gu-Skoufranis]\label{biEexpansion}Let $|c|\leq 1$ and $\mu \in\mathcal{M}_{0,c}^3(\mathbb{R}^2)$, then there exist coefficients $\{r_{k,l}\}_{k+l=3}$ such that
\[
    E_\mu(z,w)=\frac{1}{z^2}+\frac{1}{w^2}+\frac{c}{zw}+\bigg(\sum_{\substack{k+l=3\\k,l \geq 0}}\frac{r_{k,l}}{z^{k}w^{l}} + o\Big(\frac{1}{z^{k}w^{l}}\Big)\bigg)
\]
as $|z|,|w|\to \infty$ non--tangentially.
\end{proposition}
\begin{proof}
    We first prove that
    \begin{equation}\label{eq:biseries}
               G_\mu(z,w)=\frac{1}{zw}\Bigg(1+\frac{1}{z^2}+\frac{1}{w^2}+\frac{c}{zw}+\bigg(\sum_{\substack{k+l=3\\k,l \geq 0}}\frac{m_{k,l}(\mu)}{z^{k}w^{l}} + o\Big(\frac{1}{z^{k}w^{l}}\Big)\bigg)\Bigg)
    \end{equation}
    for $|z|,|w|\to \infty$
    non--tangentially. In this limit, note that there exist constants $C_1,C_2$ such that $|z-t|\geq C_1|z|, |w-s|\geq C_2|w|$ for any $s,t\in \mathbb{R}$. Using the identity
    \[
        \frac{1}{(z-t)}=\sum_{j\leq 2} \frac{t^j}{z^{j+1}}+\frac{t^{3}}{z^{3}(z-t)},
    \]
    we can therefore write
    \[
        \frac{1}{(z-t)(w-s)}=\sum_{k+\ell\leq 3} \frac{t^ks^\ell}{z^{k+1}w^{\ell+1}}+S(z,w;t,s),\quad |S(z,w;t,s)|\leq C\Big(\frac{|t|^3}{|z|^4}+\frac{|s|^3}{|w|^4}\Big)
    \]
    where $C$ only depends on $C_1,C_2$. By integrating against $\mu(s,t)$, which has finite absolute third moments, we conclude \eqref{eq:biseries}.
    The expansions in \eqref{eq:biseries} and in proposition \ref{BGtheorem} give a partial Taylor expansion of the desired order for $R_\mu(z,w)$; what is left is to argue that the coefficients match those in the proposition statement. This follows by arguing as in the proof of theorem 2.4 in \cite{Voiculescu2}, and we conclude proposition \ref{biRexpansion}. For the self-energy, note that we have expansions for $E_{\mu^{(i)}}$ and $G_{\mu^{(i)}}$, by proposition \ref{Kexpansion} and the Hamburger-Nevanlinna theorem:
    \[  
        G_{\mu^{(i)}}(z)=\sum_{k=0}^p \frac{m_k({\mu^{(i)}})}{z^{k+1}}+o\Big(\frac{1}{z^{k+1}}\Big),\quad \text{ (when ${\mu^{(i)}}$ has finite $p-$th moment)}.
    \]
    Arguing as in \cite{GuSkoufranis} (theorem 4.3) to match coefficients then yields the claim.
\end{proof}

One then defines bi--free (resp. bi--Boolean) additive convolution $\mu {\boxplus\boxplus} \nu$ (resp. $\mu{\uplus\uplus}\nu$) as the operation linearized by the two--dimensional $R$--transform (resp. self--energy), namely for which
\[
    R_{\mu{\boxplus\boxplus} \nu}(z,w)=R_\mu(z,w)+R_\nu(z,w), \quad \big(\text{resp. }
    E_{\mu{\uplus\uplus} \nu}(z,w)=E_\mu(z,w)+E_\nu(z,w)\big)
\]
where these functions are defined. Under dilation of the underlying measure, only the input of $R_\mu$ and $E_\mu$ is scaled, as opposed both the input and the function itself (e.g. in 
    \eqref{scalingRK}). 
\begin{lemma}\label{scalingRbi} Let $\mu\in\mathcal{M}(\mathbb{R}^2)$ and $\lambda \in (0,1)$. Then
\[
    R_{\lambda \mu} (z,w) = R_\mu(\lambda z,\lambda w), \quad E_{\lambda \mu}(z,w)=E_\mu( z/\lambda, w/ \lambda).
\]
\end{lemma}
\begin{proof}
    By definition of $G_\mu$, we have $G_{\lambda \mu}(z,w)=(1/\lambda)^2G_{ \mu}(z/\lambda, w/\lambda)$. It follows that
    \begin{align*}
        R_{\lambda \mu}(z,w) &= 1+\lambda z R_{{\mu}^{(1)}}(\lambda z) +\lambda wR_{{\mu}^{(2)}}(\lambda w)-\frac{(\lambda z)(\lambda w)}{G_{\mu}\big(R_{{\lambda \mu}^{(1)}}\big(\lambda z)+\frac{1}{\lambda z},R_{{\lambda\mu}^{(2)}}(\lambda w)+\frac{1}{\lambda w}\big)}
    \end{align*}
    which equals $R_\mu(\lambda z,\lambda w)$.
    The proof for $E_\mu$ is similar.
\end{proof}

\subsection{Weak convergence}

Using the bi-free R transform, weak convergence can be characterized as follows (cf. proposition 2.6 in \cite{HuangWang})
\begin{theorem}[Huang--Wang]\label{charbi}
    Let $\{\nu_n\}_{n\geq 1}$ be a sequence of Borel probability measures. Then $\nu_n$ converges weakly to a probability measure on $\mathbb{R}^2$ if and only if 
    \begin{enumerate}
        \item There exists a Stolz angle $\Delta$ such that all $R_{\nu_n}$ are defined in the product domain $\Omega=(\Delta\cup\overline{\Delta})\times (\Delta\cup\overline{\Delta})$,
        \item The pointwise limit $\lim_{n\to\infty}R_{\nu_n}(z,w)=R(z,w)$ exists for every $(z,w)$ in the domain $\Omega$, and
        \item The limit $R_{\nu_n}(-iy,-iv)\to 0$ holds uniformly in $n$ as $y,v\to 0^+$. 
    \end{enumerate}
    Moreover, if the $\nu_n$ converge weakly to $\nu$, $R_\nu=R$. 
\end{theorem}
The analogous result for the self--energy is the following (proposition 5.7 in \cite{GuSkoufranis}).
\begin{theorem}[Gu--Skoufranis]\label{charbibool} Let $\{\nu_n\}_{n\geq 1}$ be a sequence of Borel probability measures. Then the following are equivalent.
     \begin{enumerate}
         \item The sequence $\{\nu_n\}_{n\geq 1}$ converges weakly to some $\nu\in \mathcal{M}(\mathbb{R}^2)$.
         \item The pointwise limits $\lim_{n\to \infty}E_{\nu_n}(z,w)=E(z,w)$ exist for all $(z,w)\in(\mathbb{C}\setminus\mathbb{R})^2$, and the limit $E_{\nu_n}(z,w)\to 0$ holds uniformly in $n$ as $|z|,|w|\to\infty$ non--tangentially.
     \end{enumerate}
     Moreover, if these assertions hold, then $E_\nu=E$ on $(\mathbb{C}\setminus\mathbb{R})^2$.
\end{theorem}

\noindent To use theorem \ref{charbi}, we argue the existence of a common domain $\Omega_0$ on which $R_\mu$ is analytic for all $\mu\in \mathcal{M}^3(\mathbb{R}^2)$ with marginals of variance $\leq 1$ (which in particular includes $\mathcal{M}_{0,c}^3(\mathbb{R}^2)$).

\begin{proposition}\label{commondomain}
    There exists a fixed Stolz angle $\Delta$ such that $R_\mu$ is analytic in $\Omega_0 :=(\Delta\cup\overline{\Delta})\times(\Delta\cup\overline{\Delta})$ uniformly over all $\mu\in \mathcal{M}^3(\mathbb{R}^2)$ satisfying $m_2(\mu^{(i)})\leq 1, \forall i\in \{1,2\}$.
\end{proposition}
\begin{proof}
     Since the marginals $\mu^{(1)}, \mu^{(2)}$ have variance $\leq 1$, their $R$--transforms are analytic on $|z+i/4|<1/4$ and hence on $\Delta' \cup \overline{\Delta}'$ for some Stolz angle $\Delta'$ that does not depend on $\mu$.
    Now for any $r>0$, 
    \[
        \mu(\{||(t,s)||\geq r\}) \leq \mu^{(1)}\bigg(\bigg\{t:|t|\geq \frac{r}{\sqrt{2}}\bigg\}\bigg)+\mu^{(2)}\bigg(\bigg\{s:|s|\geq \frac{r}{\sqrt{2}}\bigg\}\bigg).
    \]
    It thus follows from Chebyshev's inequality and taking limits as $r\to\infty$ that the family of measures satisfying our conditions is tight. For such families of measures, $G_\mu(z,w)=(zw)^{-1}(1+o(1)) $
    uniformly in $\mu$ as $z,w\to\infty$ non--tangentially (see proposition 2.1 in \cite{HasebeHuangWang}), and we can thus pick a Stolz angle $\Delta''$ on which $G_\mu(R_{\mu^{(1)}}(z)+1/z,R_{\mu^{(2)}}(w)+1/w) $
    does not vanish for every $\mu\in \mathcal{M}_{0,c}^3(\mathbb{R}^2)$. Picking $\Delta=\Delta'\cap \Delta''$ and $\Omega_0=(\Delta\cup \overline{\Delta})\times (\Delta\cup\overline{\Delta})$ gives the result.
\end{proof}

% Given such a domain $\Omega_0$, we can always shrink it in order to ensure that $R_{\mu{\boxplus\boxplus}\nu}$ is defined $\Omega_0$ whenever $\mu, \nu \in \mathcal{M}_{0,c}^3$. In the following section (namely Proposition \ref{propertiesBF}), we assume this to be true.

\section{Extension to bi--free and bi--Boolean}
We define the \textit{standard bi--free Gaussian} $\gamma_c$ of covariance $c\in [-1,1]$  to be the measure whose $R$--transform is 
\[
    R_{\gamma_c}(z,w)=z^2+w^2+czw,
\]
which has mean $(0,0)$ and marginals of variance 1, and belongs to $\mathcal{M}_{0,c}^3(\mathbb{R}^2)$. The same is true for the \textit{standard bi--Boolean Gaussian} of covariance $c$, which is the measure $\tilde{\gamma}_c$ whose self--energy is
\[
    E_{\tilde{\gamma}_c}(z,w)=\frac{1}{z^2}+\frac{1}{w^2}+\frac{c}{zw}.
\]
When $c=0$, these measures reduce to the product measures $\rho_{sc}\otimes \rho_{sc}$ and $\rho_b\otimes\rho_b$ respectively. 

Let $\Omega_0$ be the domain from proposition \ref{commondomain}. We define the two--dimensional analogues to the distances in \eqref{metrics} on 
\begin{align*}
    d_{\text{BF}}(\mu, \nu) &= \sup_{(z,w)\in \Omega_0}\frac{|R_\mu(z,w)-R_\nu(z,w)|}{\sum_{k+l=3}|z|^k|w|^l},\quad d_{\text{BB}}(\mu, \nu) = \sup_{(z,w)\in \Omega_0}\frac{|E_\mu({1}/{z},1/w)-E_\nu(1/z,1/w)|}{\sum_{k+l=3}|z|^k|w|^l},
\end{align*}
for any $\mu,\nu \in \mathcal{M}^{3}(\mathbb{R}^2)$ with variance $\leq 1$ marginals.

\begin{theorem}[Bi-free Berry--Esseen]\label{BifBE} Let $\mu\in \mathcal{M}_{0,c}^3(\mathbb{R}^2)$. Then there exists a $C(\mu)=C>0$ such that $d_{\text{\normalfont{BF}}}(\tfrac{1}{\sqrt{n}}\mu^{{\boxplus\boxplus}n},\gamma_c)\leq {C}/{\sqrt{n}}$.
In particular, $\tfrac{1}{\sqrt{n}}\mu^{{\boxplus\boxplus}n}$ converges weakly to $\gamma_c$.
\end{theorem}
\begin{theorem}[Bi-Boolean Berry--Esseen]\label{BboolBE} Let $\mu\in \mathcal{M}_{0,c}^3(\mathbb{R}^2)$. Then there exists a $C(\mu)=C>0$ such that $d_{\text{\normalfont{BB}}}(\tfrac{1}{\sqrt{n}}\mu^{{\uplus\uplus}n},\tilde{\gamma}_c)\leq {C}/{\sqrt{n}}.$
In particular, $\tfrac{1}{\sqrt{n}}\mu^{{\uplus\uplus}n}$ converges weakly to $\tilde{\gamma}_c$.
\end{theorem}

We only prove the first of the two theorems, the proof of the second being identical. Consider the renormalization map $T^{{\boxplus\boxplus}}:\mu\mapsto \frac{1}{\sqrt{2}}(\mu{\boxplus\boxplus}\mu)$, noting that $T^{{\boxplus\boxplus}}\gamma_c = \gamma_c$, and that $\mathcal{M}_{0,c}^3(\mathbb{R}^2)$ is closed under the action of $T^{{\boxplus\boxplus}}$, as one would come to expect. By arguing as in the previous section, we straightforwardly get the following propositions.
\begin{proposition}
    $d_{\text{\normalfont{BF}}}$ is a finite metric on $\mathcal{M}_{0,c}^3(\mathbb{R}^2)$, where convergence in the metric topology implies weak convergence.
\end{proposition}
\begin{proposition}
    $T^{{\boxplus\boxplus}}$ is a contraction on $(\mathcal{M}_{0,c}^3(\mathbb{R}^2), d_{\text{\normalfont{BF}}})$ with contraction constant $1/2$.
\end{proposition}
\begin{proposition}\label{propertiesBF}
    Let $\mu, \nu, \xi, \eta \in \mathcal{M}^3(\mathbb{R}^2)$ have marginals with variance $\leq 1$, and $\lambda \in (0,1)$. Then 
    \begin{align*}
        d_{\text{\normalfont BF}}(\lambda \mu,\lambda \nu)&\leq \lambda^3d_{\text{\normalfont BF}}( \mu, \nu)
    \end{align*}
    and whenever $\mu {\boxplus\boxplus}\nu$, $\xi{\boxplus\boxplus}\eta$ have variance $\leq 1$ marginals,
    \[
         d_{\text{\normalfont BF}}(\mu {\boxplus\boxplus}\nu, \xi{\boxplus\boxplus}\eta) \leq d_{\text{\normalfont BF}}(\mu, \xi)+d_{\text{\normalfont BF}}(\nu,\eta).
    \]
\end{proposition}
The latter two follow from the proofs of \ref{contract} and \ref{properties} respectively, making the appropriate substitutions. For the first proposition, finiteness follows from the partial expansion in proposition \ref{biRexpansion}, and the fact that convergence in $d_{\mathrm{BF}}$ implies weak convergence follows from theorem \ref{charbi}.

\subsection{Other Potential Generalizations} The same proof could be adapted to any type of convolution $*$ of measures on $\mathbb{R}$ that is linearized by a set of cumulants $\{c_n\}_{n\geq 1}$, so long as the cumulants satisfy the following properties
\begin{align*}
    &c_n({\mu*\nu}) = c_n(\mu)+c_n({\nu}) &\mbox{ \text{(Additivity)}}\\
    &c_n{(\lambda \mu)} = \lambda^nc_n(\mu) &\mbox{ \text{(Homogeneity)}}
\end{align*}
and provided that one has a sufficiently well-developed analytic theory for their generating function $C(z)=\sum_{i=1}^\infty c_i z^{i-1}$. We note that these are two out of the three properties proposed by Lehner \cite{Lehner} in his axiomatization of cumulants in non--commutative probability. One can actually relax the additivity property, requiring instead that
\[  
    c_n(\mu^{*k})=kc_n(\mu)
\]
and the proof would still hold. This relaxation was used by Hasebe \cite{hasebeCumulants} to define cumulants for monotone independence, for which we cannot have additivity due to a dependence on the order of the associated random variables. 

\section*{Acknowledgments and funding}

I would like to thank my advisor, Prof.\ Jon Keating, for bringing Ott's article to my attention and for raising the question of whether the renormalization group approach could be used for the free CLT. Many thanks are also due to Dr.\ Ott himself for his valuable feedback on this work, as well as Profs.\ Roland Speicher and Michael Anshelevich for their comments (in particular for referring me to the latter's work). 
Lastly, I'm very grateful to Prof. Natasha Blitvic for introducing me to notions of independence beyond classical and free. This work was supported by the EPSRC Centre for Doctoral Training in Mathematics of Random Systems: Analysis, Modelling and Simulation (EP/S023925/1).

\nocite{*}
\bibliographystyle{abbrv}
\bibliography{biblio}

\end{document}